\documentstyle{amsart}

\newtheorem{theorem}{Theorem}
\newtheorem{lemma}{Lemma}

\newtheorem{corollary}{Corollary}
\newtheorem{definition}{Definition}
\newtheorem{example}{Example}
\newtheorem{question}{Question}

\newcommand{\om}{\omega}
\newcommand{\al}{\alpha}
\newcommand{\ah}{\aleph}
\newcommand{\emp}{\emptyset}
\newcommand{\la}{\langle}
\newcommand{\ra}{\rangle}
\newcommand{\bR}{\bf R}
\newcommand{\nio}{n \in \omega}

\newcommand{\cont}{2^{\aleph_0}}
\newcommand{\bsl}{\setminus}

\newcommand{\restrict}{|}
\newcommand{\eps}{\varepsilon}

\begin{document}

\title{Strongly almost disjoint sets and weakly uniform bases}
\author{Z. T. Balogh, S. W. Davis, 
W. Just , S. Shelah, and P. J. Szeptycki}

\thanks{ The first author's research was partially supported by NSF grant DMS-9623391.   The third author's research was done during visits at Rutgers University and The Hebrew University, Jerusalem which were supported by NSF grant DMS-9704477 and the
Landau Center. The fourth author was partially
supported by the Israel Basic Research Fund. This is publication
number~674 in
Shelah's publication list.}

\begin{abstract}
A combinatorial principle CECA is formulated and its equivalence
with GCH
+ certain weakenings of $\Box_\lambda$ for singular $\lambda$ is
proved.
CECA is used to show that  certain  ``almost point-$<
\tau$'' 
families can be refined to point-$< \tau$ families by removing a
small
set from each member of the family. This theorem in turn is used
to show
the consistency of ``every first countable $T_1$-space with a
weakly uniform
base has a point-countable base.''
\end{abstract}

\maketitle

This research was originally inspired by the following question
of Heath and 
Lindgren \cite{HL}: Does every first countable Hausdorff space
$X$
with a weakly uniform
base have a point-countable base? The answer to this question is
negative 
if MA + $\cont > \aleph_2$ is assumed (see \cite{DRW}). On the
other hand,
if CH holds and the space has at most $\aleph_\omega$ isolated
points, then the
answer is positive (see \cite{AJRS}). The starting point of this
paper was 
the observation that if in addition to GCH also the combinatorial
principle
$\Box_\lambda$ holds for every singular cardinal $\lambda$, then
no bound
on the number of isolated points is needed. An analysis of the
proof lead to
the formulation of a combinatorial principle CECA. It turns
out that CECA is  equivalent to GCH + some previously known
weakenings of 
$\Box_\lambda$, but CECA has a different flavor than
$\Box_\lambda$-principles
and may be
easier to work with. The equivalence will be shown in Section~1.

In \cite{HJS} it is shown that under certain
conditions almost disjoint families $\{ A_\alpha: \, \alpha <
\kappa\}$ can
be refined to disjoint families by removing  small sets
$A_\alpha'$ from
each $A_\alpha$. Let us say that a family 
$\{ A_\alpha: \, \alpha < \kappa\}$ is {\em point-$< \tau$} if
for every
$I \in [\kappa]^\tau$ the intersection $\bigcap_{\alpha \in I}
A_\alpha$
is empty. In particular, a family is disjoint iff it is point-$<
2$.
In Section~2 the main theorem of this paper (Theorem~\ref{main})
is derived
from CECA. Roughly speaking, Theorem~\ref{main} asserts that
certain families
$\{ A_\alpha: \, \alpha < \kappa\}$ that are ``almost point-$<
\tau$'' can
be refined to point-$< \tau$ families by removing  small sets
$A_\alpha'$ from
each $A_\alpha$.

In Section 3, some related results for almost disjoint families are proved,
and it is explored how much of Theorem~\ref{main} can be derived from GCH alone
rather than from CECA.

In Section~4 we show that the positive answer to the question of 
\cite{DRW}, and more, follows already from Theorem~\ref{main}.

\section{The Closed Continuous $\in$-chain Axiom}

\begin{definition}\label{cecadef}
{\em Let $\tau$ be a regular cardinal.
A set $M$ is {\em $\tau$-closed\/} if 
$[M]^{<\tau} \subset M$.
We say that a set $M$ is {\em weakly $\tau$-closed\/} if 
for every $I \in [M]^\tau$ there exists $J \in [I]^\tau$
such that $[J]^{<\tau} \subset M$. We say that $M$ is 
{\em weakly closed\/} if $M$ is weakly $\tau$-closed for all
regular
$\tau$ (equivalently: for all regular $\tau \leq |M|$).
An {\em $\in$-chain\/} is a
sequence $\la M_\xi: \, \xi < \alpha\ra$ such that
$M_\xi \subset M_\eta$ and $\la M_{\xi}: \, \xi \leq \eta \ra \in
M_{\eta + 1}$
for all $\xi < \eta < \alpha$.
The {\em Closed Continuous $\in$-chain Axiom for $\lambda$ and
$\tau$\/} 
(abbreviated CECA${}_\tau$($\lambda$)) is the
following statement:

{\em  For every cardinal $\Theta > \lambda$, and
for every pair of sets $A,B$ with $|A| = \lambda$, 
$|B| < \lambda$, there exists
a continuous $\in$-chain $\la M_\xi \ra_{\xi < cf (\lambda )}$ of
weakly 
$\tau$-closed
elementary submodels of $H(\Theta )$ such that\\
-- $|M_\xi | < \lambda$ for every $\xi < cf (\lambda )$;\\
-- $B \in M_0$; and\\
-- $A \subset \bigcup_{\xi < cf (\lambda )} M_\xi$.}

The {\em Closed Continuous $\in$-chain Axiom for $\lambda$\/} 
(abbreviated CECA($\lambda$)) is obtained by requiring that the
models 
$M_\xi$ in the definition of CECA${}_\tau$($\lambda$) are weakly
closed.
The {\em Closed Continuous $\in$-chain Axiom\/} 
(abbreviated CECA) asserts that CECA($\lambda$) holds for all
uncountable
cardinals $\lambda$.}
\end{definition}

\begin{lemma}\label{skolem}
Let $\lambda$ be a regular uncountable cardinal. The following are equivalent:\\
{\em (a)} $\lambda^{< \lambda} = \lambda$;\\
{\em (b)} For every set $A$ with $|A| \leq \lambda$ there exists
a weakly 
closed $M$ such that $A \subset M$\\
\hspace*{0.6cm} and $|M| = \lambda$;\\
{\em (c)} There exists a weakly $\lambda$-closed $M$ such that
$|M| = \lambda$.
\end{lemma}

\noindent
{\bf Proof:} (a) $\Rightarrow$ (b)\\

Suppose $\lambda^{< \lambda} = \lambda$, and $|A| = \lambda$.
Build recursively
an increasing sequence $\la M_\xi : \, \xi < \lambda\ra$ such
that
$A \subset M_0$, $[M_{\xi}]^{< \lambda} \subset M_{\xi + 1}$, and
$|M_\xi| = \lambda$ for all $\xi < \lambda$. 
Let $M = \bigcup_{\xi < \lambda} M_\xi$. If $\tau \leq \lambda$
is regular,
$I \in [M_\xi]^\tau$ and $K \in [I]^{< \tau}$, then $K \in
[M_\xi]^{< \tau}
\subseteq [M_{\xi}]^{< \lambda}$
for some $\xi < \lambda$, and hence $K \in M_{\xi + 1} \subseteq
M$, as
required.\\

The implications (b) $\Rightarrow$ (c) and (c) $\Rightarrow$ (a)
are obvious.\hspace{\fill}$\Box$

\begin{corollary}
{\em CECA $\Rightarrow$ GCH.}
\end{corollary}

\begin{lemma}\label{gchcase}
Assume GCH.
If $\lambda$ is an uncountable limit cardinal or if $\lambda =
\kappa^+$ for
some regular infinite cardinal $\kappa$, then CECA($\lambda$)
holds.
\end{lemma}

\noindent
{\bf Proof:} Let $\lambda$ be as above, let
$\Theta > \lambda$, and let $A,B$ be such that $|A| = \lambda$, 
$|B| < \lambda$. In both of these cases we can
construct 
a continuous $\in$-chain $\la M_\xi \ra_{\xi < cf (\lambda )}$ of

elementary submodels of $H(\Theta )$ such that\\
-- $|M_\xi | < \lambda$ for every $\xi < cf (\lambda )$;\\
-- $|M_{\xi + 1}|$ is regular and 
   $[M_{\xi + 1}]^{< |M_{\xi + 1}|}\subset M_{\xi + 1}$ 
   for every $\xi < cf (\lambda )$;\\
-- $B \in M_0$; and\\
-- $A \subset \bigcup_{\xi < cf (\lambda )} M_\xi$.

It remains to show that if $\delta < cf (\lambda )$ is limit,
then
$M_\delta$ is weakly closed. So assume $\delta < cf (\lambda )$
is limit,
let $\tau$ be a regular cardinal, and let $I \in
[M_\delta]^\tau$.

If $cf (\delta) < \tau$, then there is $\xi < \tau$ with $|M_\xi
\cap I| =
\tau$. Fix such $\xi$, and let $J = I \cap M_{\xi + 1}$. Then
$[J]^{< \tau} \subset [M_{\xi + 1}]^{< \tau} \subseteq 
[M_{\xi + 1}]^{< |M_{\xi + 1}|} \subset M_{\xi + 1} \subset
M_\delta$.
     
If $cf (\delta) \geq \tau$, then $[I]^{< \tau} \subset
\bigcup_{\xi < \delta}
M_\xi = M_\delta$, and we can take $J=I$.\hspace{\fill}$\Box$\\

Thus GCH implies that CECA($\kappa$) holds for all $\kappa$, 
except perhaps if $\kappa = \lambda^+$ for a
singular strong limit cardinal $\lambda$.  
Fortunately, it turns out that in this case 
CECA($\kappa$) is equivalent to a weakening of $\Box_\lambda$
that has been 
extensively studied by the fourth author. To prove the
equivalence,
let us introduce some notation.

We say that  $Pr(\lambda ,\tau )$ holds if 
there exists an increasing continuous chain
$\bar{N} = \la N_i : \, i < \lambda^+\ra$ such that
$|N_i| = \lambda$ for each $i < \lambda^+$,  each $N_i$ 
is weakly $\tau$-closed, and $\lambda^+ \subset \bigcup_{i <
\lambda^+} N_i$. 

We say that  $Pr'(\lambda ,\tau )$ holds if 
there exists an increasing continuous chain
$\bar{N} = \la N_i : \, i < \lambda^+\ra$ such that
$|N_i| = \lambda$ for each $i < \lambda^+$,  each $N_i$ 
is weakly $\tau$-closed with respect to sets of ordinals, 
and $\lambda^+ \subset \bigcup_{i < \lambda^+} N_i$.

Note that if $2^\lambda = \lambda^+$, then in $Pr(\lambda , \tau )$ and
$Pr' (\lambda , \tau )$ we can demand that $\bigcup_{i < \lambda^+} N_i =
H(\lambda^+ )$.

In Definition 1.4 of \cite{HJS}, 
the following principle $Sp(\sigma ,  \lambda)$ was
introduced:\footnote{Actually, the principle introduced in
\cite{HJS}
is more general and contains an extra parameter~$\tau$,
but the case $\tau = \sigma^+$ is most relevant for the results
in 
\cite{HJS} and fits most neatly into the framework of the present
paper.}
{\em There exists a sequence $\la {\cal P}_\xi : \, \xi <
\lambda^+\ra$ such that for all $\xi < \lambda^+$ we have 
${\cal P}_\xi \subset [\xi]^\sigma$ and $|{\cal P}_\xi| \leq
\lambda$;
moreover if $\xi < \lambda^+$ with $\sigma^+ = cf (\xi ) $
and $x$ is a cofinal subset of $\xi$ of cardinality $\sigma^+$, 
then $x$ can be written in the form $x = \bigcup \{ x_\nu: \, \nu
\in 
\sigma\}$ where for each $\nu \in \sigma$ 
we have $[x_\nu]^\sigma \subset \bigcup_{\eta < \xi} {\cal
P}_\eta$.}

The ideal $I[\kappa]$ was defined in \cite{108} and \cite{88a}.
We use here two equivalent definitions of $I[\kappa]$ given in 
\cite{420}.

\begin{definition}\label{ilambdadef}
{\em For a regular uncountable cardinal $\kappa$, 
let $I[\kappa ]$ be the family of all sets 
$A \subseteq \kappa$ such that the set
$\{ \delta \in A : \delta = cf (\delta )\}$  is not stationary in
$\kappa$ and
for some $\la {\cal P}_\alpha : \, \alpha < \kappa \ra$ we
have:\\
{\em (a)} ${\cal P}_\alpha$ is a family of $< \kappa$ subsets of
$\alpha$;\\
{\em (b)} for every limit 
$\alpha \in A$ such that $cf (\alpha ) < \alpha$ there
is $x \subset \alpha$ such that\\
\hspace*{0.5cm} $otp(x) < \alpha = \sup x$ and
$\forall \beta < \alpha \, (x \cap \beta \in \bigcup_{\gamma <
\alpha} {\cal P}_\gamma)$.}
\end{definition}

The following characterization of $I[\kappa]$ appears as
Claim~1.2 in 
\cite{420}. The abbreviation $nacc (C)$ stands for
``nonaccumulation points of $C$'' (in the order topology). 

\begin{lemma}\label{claim1.2}
Let $\kappa$ be a regular uncountable cardinal. Then $D \in
I[\kappa ]$
iff there exists a sequence $\la C_\beta: \, \beta < \kappa \ra$
such that:\\
{\em (a)} $C_\beta$ is a closed subset of $\beta$;\\
{\em (b)} if $\alpha \in nacc (C_\beta )$ then 
$C_\alpha = C_\beta \cap \alpha$;\\
{\em (c)} for some club $E \subseteq \kappa$, for every $\delta 
\in D \cap E$:\\
\hspace*{0.55cm}$cf (\delta) < \delta$ and $\delta =\sup
C_\delta$ and 
$otp(C_\delta) = cf (\delta)$;\\
{\em (d)} $nacc (C_\delta )$ is a set of successor ordinals.
\end{lemma} 

It is clear from the above lemma that Jensen's principle
$\Box_\lambda$
implies that $\lambda^+ \in I[\lambda^+]$.

For $\tau \leq \lambda$, let
$S^{\lambda^+}_\tau = \{ \alpha < \lambda^+: \, cf (\alpha) =
\tau\}$.

\begin{theorem}\label{shelah}
Let $\tau < \lambda$ be infinite cardinals with $\tau$ regular
and  
$\lambda$ singular strong limit. The following are 
equivalent:{\em \\
(a) $S^{\lambda^+}_\tau \in I[\lambda^+]$;\\
(b) $Pr(\lambda , \tau )$;\\
(c) $Pr'(\lambda, \tau )$;\\
(d) CECA${}_\tau$($\lambda^+$).}\\
\hspace*{0.5cm} Moreover, 
if $\tau = \sigma^+$, then each of the above is also equivalent
to:\\
{\em (e)} $Sp(\sigma , \lambda)$.
\end{theorem}

{\bf Proof:} The implications from (b) to (c) and from (d) to (b)
are
obvious.\\

(c) $\Rightarrow$ (d)\\

Let $\theta > \lambda^+$, and let $\bar{N} = \la N_i: \, i <
\lambda^+\ra$ 
exemplify $Pr'(\lambda , \tau)$. Let $\bar{M}  = \la M_\xi : \,
\xi <
\lambda^+ \ra$ 
be an increasing continuous sequence such that $\bar{N} 
\in M_0$,
and for all $\xi < \lambda^+$ we have:\\
(1) $M_\xi \prec \la H(\Theta ), <^*\ra$ (where $<^*$ is some
wellorder 
relation on $H(\Theta)$);\\
(2) $|M_\xi| = \lambda$; and\\
(3) $\bar{M}\restrict (\xi + 1) \in M_{\xi + 1}$.

Let $M =\bigcup_{\xi < \lambda^+} M_\xi$. Note that (3) implies
in particular
that $\xi \in M_{\xi + 1}$, and hence $\lambda^+ \subset M$.
Fix a bijection $h: M \rightarrow \lambda^+$ such that 
$h\restrict M_\xi \in M_{\xi+1}$ for all $\xi$. Such $h$ can be
found by
conditions (1) and (3): Recursively, let $h\restrict M_{\xi+1}$
be the
$<^*$-smallest bijection from $M_{\xi + 1}$ onto an ordinal that
extends
$h \restrict M_\xi$.
For every $\xi < \lambda^+$, let $\eta (\xi )$ be the smallest
ordinal 
$\eta \geq \xi$
such that $h[M_\eta] = \eta = M_\eta \cap \lambda^+ = N_\eta \cap
\lambda^+$. 
By (2), $\eta (\xi) < \lambda^+$ for all
$\xi$. Define:
$$E = \{ \delta \in {\bf LIM \cap \lambda^+} : \, 
 \forall \xi < \delta \, (\eta (\xi) < \delta) \}.$$
Then $E$ is a closed unbounded subset of $\lambda^+$ that
consists of fixed
points of the function $\eta$. Let
$\{ \delta_\eps: \, \eps < \lambda^+\}$ be the increasing
continuous 
enumeration of $E$. 

We show that $\la M_{\delta_\eps}: \, \eps < \lambda^+\ra$
witnesses
CECA${}_\tau$($\lambda^+$). Let $\eps < \lambda$ and let 
$I \in [M_{\delta_\eps}]^\tau$. Note that we have 
$h[I] \in [\delta_\eps]^\tau$ and $M_{\delta_\eps} \cap \lambda^+
= N_{\delta_\eps} \cap \lambda^+$. 
So there is $J \subseteq h[I]$ with $|J| = \tau$ and
$[J]^{< \tau} \subseteq N_{\delta_\eps} \cap \bf ON$. 
We distinguish two cases:\\

\noindent
{\bf Case 1:} There exists $J' \in [J]^\tau$ with $\sup J' <
\delta_\eps$.

Let $\alpha = \eta(\sup J')$. Then $\alpha < \delta_\eps$,
$M_\alpha \cap
\lambda = \alpha = h[M_\alpha]$, and 
$M_\alpha \cap \lambda^+ = N_\alpha \cap \lambda^+$.
Thus, by shrinking $J'$ if necessary, we may assume that $[J']^{<
\tau} 
\subset N_\alpha$. Let $I' = h^{-1} J'$.  
Then $I' \in [I]^\tau$. Moreover, if $K \in [I']^{< \tau}$, then 
$K = h^{-1} L$ for some $L \in [J']^{< \tau}$. Since 
$L, h\restrict M_\alpha \in M_{\alpha + 1}$
and $L \subset M_{\alpha + 1}$, it follows that $K \in M_{\alpha
+1} \subset
M_{\delta_\eps}$, as desired.\\
 
\noindent
{\bf Case 2:} There is no $J'$ as in Case 1.

Then let $I' = h^{-1} J$. Note that we must have $otp (J) = \tau$
and 
$\sup J = \delta_\eps$. Thus, if $K \in [I']^{< \tau}$, then 
$K = h^{-1} L$, where $L \in N_\alpha$ for some $\alpha <
\delta_\eps$ with
$N_\alpha = M_\alpha \cap H(\lambda^+)$. Thus, arguing as in the
previous 
case, one can show that $K \in M_{\alpha + 1}$ and hence
$[I']^{< \tau} \subset M_{\delta_\eps}$.\\

(c) $\Rightarrow$ (a)\\

Let $\bar{N} = \la N_\alpha: \, \alpha < \lambda^+\ra$ 
be a sequence that witnesses
$Pr' (\lambda , \tau)$. By thinning out the chain if necessary,
we 
may assume that $\alpha \subset N_\alpha$ for each $\alpha <
\lambda^+$.
For each $\alpha$, let ${\cal P}_\alpha =
N_\alpha \cap {\cal P } (\alpha )$. We claim that the sequence 
$\la {\cal P}_\alpha : \, \alpha < \lambda^+\ra$ witnesses that 
$S^{\lambda^+}_\tau \in I[\lambda^+]$. Condition (a) of 
Definition~\ref{ilambdadef} is obvious. To verify that (b) also
holds,
let $\alpha \in S^{\lambda^+}_\tau \bsl \{ \tau \}$. 
Pick a subset $I$ of $\alpha$ of 
order type $\tau$ such that $\sup I = \alpha$. Then $I \in
[N_\alpha \cap 
{\bf On}]^\tau$, and there exists $J \in [I]^\tau$ such that 
$[J]^{< \tau} \subset N_\alpha$. Since $\alpha$ is a limit
ordinal and 
the sequence $\bar{N}$ is continuous, the latter implies that
$[J]^{< \tau} \subset \bigcup_{\beta < \alpha} {\cal P}_\beta$.
On the other
hand, regularity of $\tau$ implies that 
$otp (J) = \tau < \alpha$ and for every $\beta < \alpha$ the set 
$J \cap \beta$ has cardinality $< \tau$. Thus condition (b) holds
and
we have shown that $S^{\lambda^+}_\tau \in I[\lambda^+]$.\\

(a) $\Rightarrow$ (d)\\

We will actually prove something slightly more general.
Assume $S^{\lambda^+}_\tau \subset D \in I[\lambda^+]$. Let 
$\la C_\alpha : \, \alpha < \lambda\ra$ and $E$ be witnesses that
$D \in 
I[\lambda^+]$ as in Lemma~\ref{claim1.2}.
Let $\mu = cf (\lambda )$, and let $(\kappa_i)_{i < \mu}$ be a
sequence of cardinals with supremum $\lambda$ and such that $\max
\{\mu , \tau\} < \kappa_0$
and $2^{\kappa_i} \leq \kappa_{i+1}$ for all $i < \mu$.
Let $A, B$ be as in the assumptions of
CECA${}_\tau$($\lambda^+$), and
let $\Theta > \lambda^+$.
Recursively construct a double sequence
$\la M^i_\alpha: \, \alpha < \lambda^+ , i<\mu\ra$ such that
conditions 
1--8 below are satisfied. For $\alpha < \lambda^+$, let 
$M_\alpha = \bigcup_{i < \mu} M^i_\alpha$. 

\begin{enumerate}
\item $M_\alpha \prec H(\Theta)$;
\item $B \subset M_0$;
\item $A \subset \bigcup_{\alpha < \lambda^+} M_\alpha$;
\item the sequence $\la M_\alpha: \, \alpha < \lambda^+\ra$ is
continuous
increasing;
\item $\la M_\beta : \, \beta \leq \alpha\ra \in M_{\alpha +1}$;
\item $|M^i_\alpha| = \kappa_i$ and $M^i_\alpha \subset
M^j_\alpha$ for $i < j$;
\item ${\cal P} (M^i_{\alpha +1}) \subset M^{i+1}_{\alpha + 1}$;
\item if $|C_\alpha| < \lambda$, then
${\cal P}(\bigcup_{\beta \in C_\alpha} M^i_\beta) \subset 
M_{\alpha + 1}$. 
\end{enumerate}

The construction is straightforward. Note that condition 8 can be
satisfied
since $|\bigcup_{\beta \in C_\alpha} M^i_\beta| = |C_\alpha|
\cdot \kappa_i$,
and $|C_\alpha| < \lambda$ for the relevant $\alpha$'s.

Now let $\la \alpha_\xi: \, \xi < \lambda^+\ra$ be the continuous
increasing
enumeration of $E$. We show that the sequence $\la
M_{\alpha_\xi}: \, \xi <
\lambda^+\ra$ witnesses CECA${}_\tau (\lambda^+)$. For this it
suffices to
verify that each $M_{\alpha_\xi}$ is weakly $\tau$-closed; the
remaining
requirements of CECA${}_\tau (\lambda^+)$ are already covered by
conditions
1--6. If $\alpha_\xi = \beta + 1$ for some~$\beta$, then
$M_{\alpha_\xi}$ is $\tau$-closed by condition~7. So consider the
case when
$\alpha_\xi$ is a limit ordinal, and let $I \in
[M_{\alpha_\xi}]^\tau$.
We distinguish four cases:\\

\noindent
{\bf Case 1:} $cf (\alpha_\xi) = \tau \neq \mu$.

Then there exist $J \in [I]^\tau$ and $i < \mu$ such that 
$J \subset \bigcup_{\beta \in C_{\alpha_\xi}} M_\beta^i$. 
We will show that $[J]^{< \tau} \subset M_{\alpha_\xi}$. Let 
$K \in [I]^{< \tau}$. Then $K \subset \bigcup_{\beta \in 
C_{\alpha_\xi} \cap \gamma} M^i_\beta$ 
for some 
$\gamma \in C_{\alpha_\xi}$. Without loss of generality, we may
assume that
$\gamma \in nacc C_{\alpha_\xi}$. Then $C_\gamma = C_{\alpha_\xi}
\cap
\gamma$.
Since $|C_{\alpha_\xi}| = \tau < \lambda$, also $|C_\gamma| <
\lambda$.
By condition~8, $K \in 
M_{\gamma + 1}$, and hence $K \in M_{\alpha_\xi}$, as required.\\

\noindent
{\bf Case 2:} $cf (\alpha_\xi) = \tau = \mu$.

We will show that $[I]^{< \tau} \subset M_{\alpha_\xi}$.
If $K \in [I]^{< \tau}$, then $K \subset M^i_{\beta + 1}$ for
some $\beta <
{\alpha_\xi}$. By condition 6, $K \in M^{i+1}_{\beta + 1} \subset
M_{\alpha_\xi}$.\\

\noindent
{\bf Case 3:} $cf(\alpha_\xi) < \tau$.

Then  there exists $\beta < \alpha_\xi$ such that 
$|I \cap M_{\beta + 1}| = \tau$. Since $M_{\beta + 1}$ is weakly
$\tau$-closed
and contained in $M_{\alpha_\xi}$, there exists $J \in 
[I \cap M_{\beta +1}]^\tau$ such that $[J]^{< \tau} \subset
M_{\alpha_\xi}$.\\ 

\noindent
{\bf Case 4:} $cf (\alpha_\xi ) > \tau$.

Let $\eta$ be the smallest limit ordinal such that $|I \cap
M_{\alpha_\eta}| =
\tau$. Then $\eta < \xi$ and $cf (\eta ) \leq \tau$. If $cf (\eta
) = \tau$, 
then we are back to Case~1 or Case~2. If $cf(\eta ) < \tau$, then
we are back
to Case~3.\\

(a) $\Rightarrow$ (e)\\

Assume that $\tau = \sigma^+$, and let $\la M_\alpha^i: \, \alpha
< \lambda^+,
i < \mu\ra$ be as in the previous part of the proof, i.e.,
such that conditions 1--8 hold, and let the sequence $\la
\alpha_\xi: \, 
\xi < \lambda^+\ra$ be defined as above. 
For each $\xi < \lambda^+$, let
${\cal P}_\xi = [\xi]^\sigma \cap M_{\xi + 1}$. 

Now suppose that $cf (\xi ) = \tau$, and let $x \subset \xi$ be
cofinal
of cardinality $\tau$. We distinguish two cases:\\

\noindent
{\bf Case 1:} $\sigma < \mu$.

Then we can let $x_\nu = x$ for all $\nu < \sigma$. Each $y \in
[x]^\sigma$
is contained in $M^i_{\alpha_\eta + 1}$ for some $i < \mu$ and
$\eta < \xi$.
Thus $[x]^\sigma \subset \bigcup_{\eta < \xi} {\cal P}_\eta$, as
required.\\

\noindent
{\bf Case 2:} $\sigma \geq \mu$.

Then let $x_\nu = x \cap \bigcup_{\beta \in C_{\alpha_\xi}}
M_\beta^\nu$
for $\nu < \mu$ and $x_\nu = \emp$ for $\nu \geq \mu$. If $y \in 
[x_\nu]^\sigma$, then $y \subset \gamma \cap 
\bigcup_{\gamma \cap C_{\alpha_\xi}}$ for some $\gamma \in nacc 
C_{\alpha_\xi}$, and hence $y \in {\cal P}_\eta$ for some
$\eta < \xi$ with $\gamma < \alpha_\eta$.\\

(e) $\Rightarrow$ (a)\\
 
Assume $\tau = \sigma^+$, and let $\la {\cal P}_\xi: \, \xi <
\lambda^+
\ra$ witness that $Sp(\sigma, \lambda)$ holds. Let
$\xi \in S^{\lambda^+}_\tau \bsl \{ \tau\}$. Now if $x \subset
\xi$
is any cofinal subset of $\xi$ of order type $\tau$, 
then $x = \bigcup_{\nu \in \sigma} 
x_\nu$, where $[x_\nu]^\sigma \subset \bigcup_{\eta < \xi} 
{\cal P}_\eta$ for each $\nu < \sigma$. At least one of these
$x_\nu$'s must
be cofinal in $\xi$ and of order type $\tau$, 
and this $x_\nu$ is exactly as required in 
Definition~\ref{ilambdadef}.\hspace{\fill}$\Box$\\

Note that
$\forall \tau \in {\bf Reg} $ CECA${}_\tau (\kappa )$ does not
always imply 
CECA($\kappa$). For example, if $2^{\aleph_n} = \aleph_{n+2}$ for
all
$\nio$, then there are no weakly closed models of cardinality
$\aleph_n$
for $n \in (0, \omega)$, and thus CECA($\aleph_\omega$) fails.
But if $\tau = \aleph_n$, then for each $m \geq n + 1$ there are
plenty
of $\aleph_{n}$-closed models of cardinality $\aleph_m$, and thus
CECA${}_{\aleph_n}$($\aleph_\omega$) holds.
This anomaly cannot happen if $\kappa$ is a successor cardinal:
If $\kappa = \lambda^+$ for a regular $\lambda$, then 
CECA${}_\lambda$($\kappa$) implies the existence of weakly
$\lambda$-closed
models of cardinality $\lambda$, which in turn implies 
that $\lambda^{<\lambda}$. Now the proof of Lemma~\ref{gchcase}
can be adapted
to derive CECA($\kappa$). If $\kappa$ is the successor of a
singular limit
cardinal $\lambda$ and $\tau < \lambda$, then 
CECA${}_{\tau^+}$($\kappa$) implies that $2^\tau \leq \lambda$. 
Since for cofinally many $\tau$ we also will have
$\tau > cf (\lambda)$, K\"onig's Theorem implies that $2^\tau$ is
strictly
less than $\lambda$. In other words,
if CECA${}_\tau$($\lambda^+$) holds and $\lambda$ is singular,
then~$\lambda$ must be a strong limit cardinal.
For such $\lambda$ we have the following  corollary to the proof
of 
Theorem~\ref{shelah}.

\begin{corollary}\label{equiv}
Let $\lambda$ be a singular strong limit cardinal. Then the
following
are equivalent:\\
{\em (a)} $\lambda^+ \in I[\lambda^+]$;\\
{\em (b)} $S^{\lambda^+}_\tau \in I[\lambda^+]$ for all regular 
$\tau < \lambda$;\\
{\em (c) CECA($\lambda^+$)};\\
{\em (d) CECA${}_\tau (\lambda^+)$} holds for all regular 
$\tau \in \lambda$.
\end{corollary}

{\bf Proof:} The equivalence between (a) and (b) follows from
Shelah's 
observation that $I[\lambda^+]$ is a normal, and hence
$\lambda^+$-complete
ideal (see \cite{108}, \cite{88a}, or \cite{420}).
The equivalence between (b) and (d) was established in
Theorem~\ref{main}.
The implication (c) $\Rightarrow$ (d) is obvious. To see that 
(a) implies (c), note that the last part of the proof of
Theorem~\ref{main}
shows that if $\lambda^+ \in I[\lambda^+]$,
then the sequence $\la M_{\alpha_\xi}: \, \xi < \lambda^+\ra$
constructed from
$D = \lambda^+$ witnesses 
CECA${}_\tau (\lambda^+)$ simultaneously for all $\tau \in
\lambda \cap
\bf Reg$.\hspace{\fill}$\Box$

\section{The main theorem}

Let $\tau$ be a cardinal. Recall that a sequence of sets 
$\la  A_\alpha \ra_{\alpha < \kappa}$ is 
{\em point $< \tau$} if for every $I \in [\kappa]^\tau$ 
the intersection $\bigcap_{\alpha \in I} A_\alpha$
is empty. 

\begin{theorem}\label{main} Assume CECA. 
Suppose that $\sigma, \tau$ are regular infinite cardinals, and
let
$\la A_\alpha\ra_{\alpha < \kappa}$ be a sequence of
(not necessarily distinct) sets such that one of the following
conditions
is satisfied:\\
$(1.1) \mbox{ For every } I \in [\kappa]^\tau
\mbox{ there is }
J \in [I]^{< \tau} \mbox{ such that } |\bigcap_{\alpha \in J}
A_\alpha| <
\sigma.$\\
$(1.2) \mbox{ For every } I \in
[\kappa]^\tau
\exists J\in [I]^{< \tau} \, |\bigcap_{\alpha \in J} A_\alpha |
\leq \sigma 
\mbox{ and }$\\
\hspace*{0.9cm}$ \forall S \in [\bigcup_{\alpha < \kappa} A_\alpha]^\sigma \,
(| \{ \alpha : \, S \subset A_\alpha \} | \leq \sigma ).$

Then there exist $\la A_\alpha'\ra_{\alpha < \kappa}$ such that
$|A_\alpha'| \leq \sigma$ for each $\alpha < \kappa$ and
the sequence $\la A_\alpha \bsl A_\alpha'\ra_{\alpha < \kappa}$
is 
point-$< \tau$.
\end{theorem}

{\bf Proof:} 
Suppose the theorem is false,  let $\kappa$ be the smallest
cardinal for which the theorem fails,  fix a counterexample
$\bar{A}= 
\la A_\alpha\ra_{\alpha < \kappa}$ and $\sigma , \tau$ that
witness this
fact. Throughout the proof, let $\Theta$ denote a ``sufficiently
large''
cardinal, and let $A = \bigcup_{\alpha < \kappa} A_\alpha$. 
In the proof we will consider all possible ways in which 
$\kappa ,\sigma, \tau$ can be related to each other, and we will
derive a
contradiction in each case. To begin with, note that we may
without loss of
generality assume that $\tau \leq \kappa$; otherwise the
conclusion of the
theorem is vacuously true. Now let us eliminate the case $\sigma
\geq \kappa$.

\begin{lemma}\label{sigmalemma}
Suppose that $M \prec H(\Theta )$ is weakly $\tau$-closed with respect to
sets of ordinals, i.e.,
for every $I \in [M\cap \kappa]^\tau$ there exists $J \in
[I]^\tau$ with
$[J]^{< \tau} \subset M$. Moreover, suppose $\bar{A} \in M$ and
$\sigma \subset M$. Then the sequence $\la A_\alpha \bsl
M\ra_{\alpha
\in \kappa \cap M}$ is point-$< \tau$. 
\end{lemma}

\begin{corollary}\label{sigmacor}
$\sigma < \kappa$.
\end{corollary}

{\bf Proof:} Let $M \prec H(\Theta )$ be weakly closed and
such that $\sigma \subset M$,
$|M| = \sigma$, and $M$ contains everything relevant.
For each $\alpha < \kappa$, let $A_\alpha' = M$. If $\sigma \geq
\kappa$, then
$\kappa \cap M = \kappa$, and 
 Lemma~\ref{sigmalemma} implies the conclusion of 
Theorem~\ref{main}.\hspace{\fill}$\Box$\\

{\bf Proof of Lemma~\ref{sigmalemma}:} Suppose towards a
contradiction that
there are $a$ and $I \in [\kappa \cap M]^\tau$ such that 
$a \in \bigcap_{\alpha \in I} A_\alpha\bsl M$. Since $M$ is
weakly $\tau$-closed with respect to sets of ordinals,\footnote{If $M$ has 
only the property that for every subset $I \in [M \cap \kappa]^{\tau^+}$ there
exists $J \in [I]^\tau$ with $[J]^{< \tau} \subset M$, then the argument
presented here shows that the sequence 
$\la A_\alpha \bsl M\ra_{\alpha \in \kappa \cap M}$ is point-$< \tau^+$.}
by passing to a subset of $I$ if necessary, we may assume that 
$[I]^{< \tau} \subset M$. Since (1.1) or (1.2) hold, 
we can find $J \in [I]^{<
\tau}$ such that 
$$(2)\qquad B = |\bigcap_{\alpha \in J} A_\alpha| \leq \sigma.$$
Since $\bar{A}, J \in M$, it follows that $B \in M$. Since
$\sigma \subset M$
and $|B| \subset \sigma$, we conclude that $B \subset M$, which
contradicts
the assumption that $a \in B \bsl M$.\hspace{\fill}$\Box$\\

Now let us reveal how CECA will be used in the remainder of this
proof.

\begin{lemma}\label{chainlemma}
Suppose that there is a continuous increasing chain $\la
M_\xi\ra_{\xi < 
cf (\kappa )}$ of weakly closed elementary submodels of $H(\Theta
)$ 
such that $|M_\xi| < \kappa$ for each $\xi < cf (\kappa)$,
that $M_0$ contains everything relevant, $\sigma \subset M_0$, 
$[\kappa]^{<cf(\kappa )} \subset \bigcup_{\xi < cf (\kappa )}
M_\xi$, 
and for every 
$\xi < \kappa$:\\
{\em (i)} The sequence 
$\la A_\alpha \bsl M_\xi\ra_{\alpha \in M_\xi \cap \kappa}$
is point-$< \tau$;\\
{\em (ii)} $\alpha \in \kappa \bsl M_\xi$ implies $|M_\xi \cap
A_\alpha| \leq 
\sigma$;\\
{\em (iii)} If $\alpha \in M_\xi$ and $|A_\alpha | \leq |M_\xi|$,
then
$A_\alpha \subset M_\xi$.\\
Then there are $A_\alpha' \in [A_\alpha]^{\leq \sigma}$ for
$\alpha < \kappa$
such that the sequence $\la A_\alpha \bsl A_\alpha'\ra_{\alpha ,
\kappa}$
is point-$< \tau$.
\end{lemma}

{\bf Proof:} Let the $M_\xi$'s and $A_\alpha$'s be as in the
assumption.
By the choice of $\kappa$ and (ii) we can pick, for every $\xi <
\kappa$ and
$\alpha \in \kappa \cap M_{\xi + 1} \bsl M_\xi$, an $A_\alpha'
\in 
[A_\alpha]^{\leq \sigma}$ in such a way that $A_\alpha \cap M_\xi
\subset
A_\alpha'$ and the sequence $\la A_\alpha \bsl
A_\alpha'\ra_{\alpha \in \kappa
\cap (M_{\xi + 1} \bsl M_\xi)}$ is point-$< \tau$.

We claim that the sequence $\la A_\alpha \bsl
A_\alpha'\ra_{\alpha < \kappa}$
is point-$<\tau$. To see this, let $a \in A$, 
and let $M^* =\bigcup_{\xi < cf (\kappa )} M_\xi$.\\

{\bf Case 1:} $a \in  M^*$.

Then there exists a unique $\xi < cf (\kappa )$ such that $a \in
M_{\xi + 1}
\bsl M_{\xi}$. For this $\xi$, the following hold:
\begin{itemize}
\item By (i), 
$|\{ \alpha\in \kappa \cap M_\xi : \, a \in A_\alpha \} | <
\tau$.
\item By construction, $|\{ \alpha \in \kappa \cap (M_{\xi + 1}
\bsl M_\xi ) : \,
   a \in A_\alpha \bsl A_\alpha'\} | < \tau$.
\item If $\alpha \in \kappa \bsl M_{\xi + 1}$, then $A_\alpha
\cap M_{\xi + 1} 
   \subset A_\alpha'$, and hence no $\alpha \in \kappa \bsl
M_{\xi + 1}$ 
   satisfies $a \in A_\alpha \bsl A_\alpha'$.
\end{itemize}

It follows that the sequence $\la A_\alpha \bsl
A_\alpha'\ra_{\alpha < \kappa}$
is point-$< \tau$ at $a$.\\

{\bf Case 2:} $a \notin M^*$.

We show that $M^*$ is weakly $\tau$-closed with respect to sets of ordinals.
Let $I \in [M^* \cap \kappa]^\tau$. If $\tau = cf (\kappa )$,
then 
$[\kappa]^{< \tau } \subset M^*$ and hence $[I]^{< \tau} \subset
M^*$. 
If $\tau \neq cf (\kappa )$, then there exists $\xi < cf (\kappa
)$ with
$|M_\xi \cap I| = \tau$, and $J = M_\xi \cap I$ is as required.
Now Lemma~\ref{sigmalemma} implies that the sequence
$\la A_\alpha \bsl A_\alpha'\ra_{\alpha < \kappa}$ is point-$<
\tau$ at 
$a$.\hspace{\fill}$\Box$\\

Since we may assume $\sigma < \kappa$, CECA clearly implies the
existence
of a sequence $\la M_\xi\ra_{\xi < cf (\kappa )}$ that satisfies
all 
unnumbered assumptions of Lemma~\ref{chainlemma}. By
Lemma~\ref{sigmalemma},
this sequence will also satisfy condition (i). How can we make
sure that
condition (ii) also holds? This depends on $\tau$.

\begin{lemma}\label{taulemma}
Suppose that $M \prec H(\Theta )$ is weakly closed, $\bar{A} \in
M$, and
$\max \{ \sigma , \tau\} \subset M$. Then $\alpha \in \kappa \bsl
M$ implies 
$|A_\alpha \cap M| \leq \sigma$.
\end{lemma}

{\bf Proof:} Assume towards a contradiction that $\alpha \in
\kappa \bsl M$
is such that $|A_\alpha \cap M| \geq \sigma^+$. Since $M$ is
weakly closed,\footnote{Note that in this argument, as well as in the proof 
of the next lemma, only the following consequence of weak closedness is used:
If $I \in [M]^{\sigma^+}$, then there is $J \in [I]^\sigma$ such that 
$J \in M$.}
there exists $Z \in [A_\alpha \cap M]^\sigma \cap M$. We are
going to prove
that the set $I = \{ \beta \in \kappa: \, Z \subset A_\beta\}$
has cardinality
$> \max \{ \sigma ,\tau\}$, 
in contradiction with conditions (1.1) and (1.2) of
Theorem~\ref{main}. 
Suppose
that $|I| \leq \max \{ \sigma ,\tau\}$. 
Note that since $\kappa, \bar{A}, Z \in M$, it follows
that $I \in M$, and there is a one-to-one function $f \in M$ that
maps
$I$ into $\max \{ \sigma , \tau\} \subset M$. It follows that $I
\subset M$.
On the other hand, $\alpha \in I \bsl M$, which gives a 
contradiction.\hspace{\fill}$\Box$\\

It follows from Lemma~\ref{taulemma} that if $\tau < \kappa$,
then we can find
a sequence $\la M_\xi\ra_{\xi < cf (\kappa )}$ that satisfies the
assumption
of Lemma~\ref{chainlemma}. So it remains to prove
Theorem~\ref{main} for the 
case $\sigma < \tau = \kappa$. If $\tau$ happens to be equal to
$\sigma^+$,
then any chain we get from CECA will satisfy $|M_\xi| \leq
\sigma$, and we
get condition (ii) for free. If $\sigma^+ < \tau = \kappa$, then
we need to
take advantage of GCH. The following lemma shows how to handle
this last
remaining case. 

\begin{lemma}\label{gchlemma}
Suppose $\kappa = \tau \geq \sigma^{++}$, and let 
$\la M_\xi\ra_{\xi < \kappa}$ be a continuous increasing
sequence of elementary
submodels of $H(\Theta )$ such that each $M_\xi$ has cardinality
less 
than $\kappa$ and $[\kappa]^{< \kappa} \subset \bigcup_{\xi <
\kappa} M_\xi$.
Then there exists a continuous subsequence $\la
M_{\xi_\nu}\ra_{\nu < \kappa}$
such that
\begin{enumerate}
\item[(3)] $\nu < \kappa$ and $\alpha \in \kappa \bsl
M_{\xi_\nu}$ implies $|A_\alpha \cap M_{\xi_\nu} | \leq \sigma .$\end{enumerate}
\end{lemma}

{\bf Proof:} Let $A = \bigcup_{\alpha < \kappa} A_\alpha$. Note
that since
$\kappa \geq \tau$, conditions (1) in Theorem~\ref{main} implies
that for
every $S \in [A]^\sigma$ the set 
$h(S) = \{ \alpha < \kappa : \, S \subset A_\alpha \}$ has
cardinality
$< \kappa$. Since $\kappa \geq \sigma^{++}$ and $\kappa\ (=
\tau)$  is 
regular, GCH implies that there is a continuous subsequence
$\la M_{\xi_\nu}\ra_{\nu < \kappa}$ of $\la M_\xi\ra_{\xi <
\kappa}$ such
that
\begin{enumerate}
\item[(4)] For every $\nu < \kappa$ and every 
$S \in [A]^\sigma \cap M_{\xi_{\nu + 1}}$ we have $h(S) \subset M_{\xi_{\nu + 1}}.$
\end{enumerate}
We claim that $\la M_{\xi_\nu}\ra_{0 < \nu < \kappa}$ is as
required.
Let $N_\nu$ denote $M_{\xi_\nu}$. Suppose towards a contradiction
that
$0 < \mu < \kappa$ and $\alpha \in \kappa \bsl N_\mu$ are such
that 
$|A_\alpha \cap N_\mu| \geq \sigma^+$. We distinguish two
cases.\\

{\bf Case 1:} $\mu = \nu + 1$ for some $\nu$.

Since $N_{\nu + 1}$ is weakly closed, there exists $S \in
[A_\alpha]^\sigma
\cap N_{\nu + 1}$. By (4), $h(S) \subset N_{\nu + 1}$. Since
$\alpha$ in
$h(S)$, we have $\alpha \in N_{\nu + 1} = N_\mu$, contradicting
the choice
of $\alpha$.\\

{\bf Case 2:} $\mu$ is a limit ordinal.

We will show that 
\begin{enumerate}
\item[(5)] There are $\nu < \mu$  and $S \in [A_\alpha]^\sigma$ such that $S\in N_{\nu +
1}.$
\end{enumerate}
Indeed, since $N_\nu$ is weakly closed and $|A_\alpha \cap N_\nu
| \geq
\sigma^+$, we can pick $S \in [A_\alpha ]^\sigma \cap N_\nu$.
Since
$N_\nu = \bigcup_{\nu < \mu} N_\nu$, there is $\nu < \mu $ such
that 
$S \in N_{\nu + 1}$. Now we can derive a contradiction as in the
previous
case.\hspace{\fill}$\Box$ $\Box$

\section{Related Theorems and Examples}

We first present a modification of Theorem~\ref{main} exhibiting
when we can obtain point-$<\tau$ families where $\tau$ is finite.
Note that a family of sets is point-$<2$ if and only if it is
pairwise disjoint.

\begin{theorem}\label{main2} 
Assume CECA and let $n \geq 1$.
Suppose that $\sigma$ is a regular cardinal and let
$\la A_\alpha\ra_{\alpha < \kappa}$ be a sequence of
(not necessarily distinct) sets such that $|A_\alpha|\leq
\sigma^{+n}$ and one of the following conditions
is satisfied:\\
$(1.3)\  |A_\alpha \cap A_\beta| < \sigma
\mbox{ \ for all \ } \alpha < \beta < \kappa .$\\
$(1.4) \  |A_\alpha \cap A_\beta | \leq \sigma  \mbox{ \ for all
\  }
\alpha < \beta < \kappa\, \mbox{ and } $\\
\hspace*{0.9cm}$\forall S \in [\bigcup_{\alpha < \kappa} A_\alpha
]^\sigma \, 
(|\{ \alpha : \, S \subset A_\alpha \} | \leq \sigma ).$

Then there exist $\la A_\alpha'\ra_{\alpha < \kappa}$ such that
$|A_\alpha'| \leq \sigma$ for each $\alpha < \kappa$ and
the sequence $\la A_\alpha \bsl A_\alpha'\ra_{\alpha < \kappa}$
is 
point-$< n+1$.
\end{theorem}

{\bf Proof:} First apply Theorem~\ref{main} for $\tau=\aleph_0$
to get sets $A_{\alpha}''\in [A_{\alpha}]^\sigma$ such that $\la
A_\alpha \bsl A_\alpha''\ra_{\alpha < \kappa}$ is 
point-finite. The step from point-finite to point-$<n+1$ can be made
in ZFC:

\begin{lemma}\label{lemmamain} Let $n \geq 1$.
Suppose that $\sigma$ is a regular
infinite cardinal and that $\{A_{\alpha}:\alpha<\kappa\}$ is a
point-finite family of sets satisfying $|A_{\alpha}|\leq
\sigma^{+n}$ for each $\alpha$ and $|A_{\alpha}\cap A_{\beta}|\leq
\sigma$ for each pair of distinct $\alpha,\beta\in\kappa$. Then
there are $A_{\alpha}'\in[A_{\alpha }]^{\sigma }$ for each
$\alpha $ such that $\{A_{\alpha}\setminus
A_{\alpha}':\alpha<\kappa\}$ is point-$< n+1$.
\end{lemma}

{\bf Proof Lemma~\ref{lemmamain}:} First let us prove the 
lemma for the case $n = 1$. Let $\kappa$ be the minimal cardinal
such that
there is a family of size $\kappa$ that forms a counterexample to
the theorem. Clearly $\kappa >\sigma $. First assume that
$\kappa>\sigma ^{+}$. Fix an $\in$-chain of elementary submodels
$\{M_{\xi}:\xi<cf(\kappa)\}$ such that 
\begin{enumerate}
\item $\sigma ^{+}\subseteq M_0$ and $|M_{\xi}| < \kappa$, for
each $\xi<\kappa$.
\item $\{A_{\alpha}:\alpha<\kappa\}\in M_0$ and 
\item $\{A_{\alpha}:\alpha<\kappa\}\subseteq
\bigcup\{M_{\xi}:\xi<cf(\kappa)\}$. 
\end{enumerate}
Using the minimality of $\kappa$, for each $\xi\in cf(\kappa)$
fix a sequence $\{A_{\alpha}':\alpha\in \kappa\cap
(M_{\xi+1}\setminus M_{\xi})\}$ of sets of size $\sigma$ such that
$\{A_{\alpha}\setminus A_{\alpha}':\alpha\in \kappa\cap
M_{\xi+1}\}$ is pairwise disjoint. Note that by point-finiteness
if $\alpha\not\in M_{\xi}$ then $A_{\alpha}\cap
M_{\xi}=\emptyset$. Also $\sigma ^{+}\subseteq M_{\xi}$ implies
that $A_{\alpha}\subseteq M_{\xi}$ for each $\alpha\in M_{\xi}$.
Therefore $\{A_{\alpha}\setminus A_{\alpha}':\alpha\in \kappa\}$
is pairwise disjoint. 

The case for $\kappa=\sigma^+ $ is similar. Fix an $\in$-chain
as above of elementary submodels $\{M_{\xi}:\xi<cf(\kappa)\}$ of
size $\sigma $ so that $\sigma\subseteq M_{0}$ and let
$A_{\alpha}'=A_{\alpha}\cap M_{\xi+1}$ where $\alpha\in
M_{\xi+1}\setminus M_{\xi}$. To see that $\{A_{\alpha}\setminus
A_{\alpha}':\alpha\in \kappa \}$ is pairwise disjoint it suffices
to observe that, as above, if $\alpha\not\in M_{\xi}$ then
$A_{\alpha}\cap M_{\xi}\not=\emptyset$ and, in addition, if
$\alpha,\beta\in M_{\xi}$ then $A_{\alpha}\cap A_{\beta}\subseteq
M_{\xi}$. 

Now suppose $m > 1$ and the lemma is true for all
$n < m$. Let $\la A_\alpha \ra_{\alpha < \kappa}$ be a sequence
of sets such that $|A_\alpha| \leq \sigma^{+m}$ for each $\alpha$ and
$|A_\alpha \cap A_\beta| \leq \sigma$ for each pair of distinct $\alpha, 
\beta \in \kappa$. Let $\nu = \sigma^{+(m-1)}$. Then 
$\sigma^{+m} = \nu^+$, and hence the lemma for $n = 1$
implies that there are $A_\alpha'' \in [A_\alpha]^\nu$ for each $\alpha$
such that the family $\{A_\alpha \bsl A_\alpha' : \, \alpha < \kappa\}$ is
pairwise disjoint. Since the lemma is true for $n = m-1$, 
there are $A_\alpha' \in [A_\alpha'']^\sigma$ such that the family
$\{ A_\alpha'' \bsl A_\alpha': \, \alpha < \kappa\}$ is point-$< m$. Now we
claim that the family $\{A_\alpha \bsl A_\alpha': \, \alpha < \kappa\}$ is
point-$< m+1$. To see this, suppose $a \in ((A_{\alpha_0}\bsl A_{\alpha_0}) 
\cap \dots \cap (A_{\alpha_k}\bsl A_{\alpha_k}))$. Then there are at most
$m-1$ indices $\alpha_i$ such that $a \in A_{\alpha_i}'' \bsl A_{\alpha_i}'$,
since the family $\{ A_\alpha'' \bsl A_\alpha': \, \alpha < \kappa\}$ is 
point-$< m$. Similarly, since the family
$\{ A_\alpha \bsl A_\alpha'': \, \alpha < \kappa\}$ is disjoint, there is at
most one $\alpha_i$ with $a \in A_{\alpha_i} \bsl A_{\alpha_i}''$. Thus, $a$
is a member of $A_\alpha \bsl A_\alpha'$ for at most $m$ indices $\alpha$,
and we are done. 
\hfill$\Box$ $\Box$\\

The following example shows that Lemma~\ref{lemmamain} is in a sense the best
possible.

\begin{theorem}\label{paulexample}
For each $n\in\omega$ there is a point-$<n+1$ family
${\cal A}_n$ of sets of size
$\aleph_n$ such that 
\begin{enumerate}
\item[(a)] $A\cap B$ is finite for each pair of distinct 
$A,B\in {\cal A}_n$.
\item[(b)]  ${\cal A}_n$ has no point-$<n$ refinement of the form $
\{A\setminus A':A\in {\cal A}_{n}\}$ where each $A'$ is
countable.
\end{enumerate}
\end{theorem}

{\bf Proof:} For each $n\in \omega $, each $i\in n$ and each
$s\in \omega _{n}^{n\setminus\{i\}}$ let $A_{s}=\{t\in \omega_{n}^{n}:\,
t|(n\setminus\{i\})=s\}$. Let ${\cal A}_{n}= \bigcup_{i \in n} \{ A_{s}:s\in
\omega _{n}^{n\setminus\{i\}}\}$. We will show that ${\cal A}_{n}$
is the required family for each $n\in\omega $.

($n=2$): ${\cal A}_{2}$ is the family of rows and columns in
$\omega_2\times\omega_2$. Clearly this family is point-$<3$ and
the intersection of any two elements is either disjoint or
contains one element. To see that (b) holds,
suppose that $A'\in [A]^{\omega }$ for each $A\in {\cal A}_{2}$.
Let $A_{0,\alpha }=\{(\alpha,\beta):\beta\in \omega_{2}\}$ and let
$A_{1,\alpha }=\{(\beta,\alpha):\beta\in \omega _{2}\}$.  Fix an
uncountable $\alpha<\omega _{2} $ so that $A_{0,\beta}'\subseteq
\alpha ^{2}$ for each $\beta \in \alpha $.  Consider $A_{1,\alpha
}$. As $\alpha $ is uncountable, we can fix $\beta <\alpha $ so
that $(\beta ,\alpha )\not\in A_{1,\alpha} '$. But then also
$(\beta ,\alpha )\not\in A_{0,\beta }'$ so the refined family is
not point-$<2$.  

Assume by induction that $n> 2$ and ${\cal A}_{k}$ is as
required for each $k<n$. Clearly ${\cal A}_{n}$ is point-$<n+1$
and almost disjoint. So suppose that $A'\in [A]^{\omega }$ for
each $A\in {\cal A}_{n}$. Fix $\alpha <\omega _{n}$ such that
$|\alpha |=\omega _{n-1}$ and so that $A_{s}'\subseteq \alpha
^{n}$ for each $i\in n$ and each $s\in \alpha^{n\setminus
\{i\}}$. Let 
$$S=\{t\in \omega _{n}^{n}:\, t(n-1)=\alpha \mbox{ and
}t|(n-1)\in \alpha ^{n-1}\}.$$ 
For each $i<n-1$, and each $s\in
\alpha ^{n\setminus\{i\}}$ such that $s(n-1)=\alpha $, let
$B_{s}=\{t\in S:t|(n-1\setminus\{i\})=s\}\subseteq A_{s}$. By the
induction hypothesis, $\{B_{s}\setminus A_{s}':
s\in\alpha^{n\setminus\{i\}},i\in n-1 \mbox{ and }s(n-1)=\alpha\}$
is not point-$<n-1$. Therefore, there is $t\in S$ such that for
each $i<n-1$, $t\in A_{t|(n\setminus\{i\})}\setminus
A_{t|(n\setminus\{i\})}' $. Finally, since $t|(n-1)\in \alpha ^{n-1}$
we have $A_{t|n-1}'\subseteq \alpha ^{n}$. Therefore $t\in A_{t|(n-1)}
\bsl A_{t \restrict (n-1)}'$, so $\{
A\setminus A':A\in {\cal A}_{n}\}$ is not point-$<n$.
\hfill$\Box$\\

By taking ${\cal A}=\bigcup_{n\in \omega }{\cal A}_{n}$ we obtain
a family satisfying the following:

\begin{corollary} There is a point-finite, almost disjoint family
${\cal A}$ such that for each sequence $\{A':A\in {\cal A}\}$ of
countable sets, the family $\{A\setminus A':A\in{\cal A}\}
$ is not point-$<n$ for any $n\in \omega$.\hspace{\fill}$\Box$
\end{corollary}

We now present the main applications of Theorem~\ref{main} and
Theorem~\ref{main2}. Letting $\sigma = \tau = \aleph_0$ in (1.1),
we get:

\begin{corollary}\label{cor1}
Assume CECA. Suppose that $\la A_\alpha\ra_{\alpha < \kappa} $ is
a sequence 
of sets such that every $I \in [\kappa ]^{\aleph_0}$ has a finite
subset 
$J \subset I$ such that $|\bigcap_{\alpha \in J} A_\alpha| <
\aleph_0$. Then
there are $A_\alpha' \in [A_\alpha]^{\leq \aleph_0}$ for $\alpha
< \kappa$
such that the sequence $\la A_\alpha \bsl A_\alpha'\ra_{\alpha <
\kappa}$
is point-finite.
\end{corollary}

Letting $\sigma = \aleph_0$ and $\tau = \aleph_1$ in (1.1), we
get:

\begin{corollary}\label{cor2}
Assume CECA. Suppose that $\la A_\alpha\ra_{\alpha < \kappa} $ is
a sequence 
of sets such that or every countably infinite set $S$,
the set $\{ \alpha < \kappa : \, S \subset A_\alpha \} $ is
countable.
Then there are $A_\alpha' \in [A_\alpha]^{\leq \aleph_0}$ for
$\alpha < \kappa$
such that the sequence $\la A_\alpha \bsl A_\alpha'\ra_{\alpha <
\kappa}$
is point-countable.
\end{corollary}

Let us restate the first part of Theorem~\ref{main2} in a more conventional
way:

\begin{corollary}\label{cor3}
Assume CECA. Let $\sigma$ be a regular infinite cardinal.
Suppose that $\la A_\alpha\ra_{\alpha < \kappa} $ is a sequence 
of sets of cardinality $\sigma^+$ each
such that $A_\alpha \cap A_\beta$ has cardinality less than 
$\sigma$ for all $\alpha < \beta < \kappa$. 
Then there are $A_\alpha' \in [A_\alpha]^{\leq \sigma}$ for
$\alpha < \kappa$
such that the sequence $\la A_\alpha \bsl A_\alpha'\ra_{\alpha <
\kappa}$
consists of pairwise disjoint sets.
\end{corollary}

The consistency of Corollary~\ref{cor3} 
 was already derived from a statement
similar to CECA as Theorem~2.6 in \cite{HJS}. 

Condition (1.4) allows us in certain circumstances to
relax the assumption that the intersection of each two sets
$A_\alpha , 
A_\beta$ has cardinality strictly less than $\sigma$. In
particular, if
$\sigma = \aleph_0$, then we get:

\begin{corollary}\label{cor4}
Assume CECA. Suppose that $\la A_\alpha\ra_{\alpha < \kappa} $ is
a sequence 
of sets of cardinality $\aleph_1$ 
such that for every countably infinite set $S$,
the set $\{ \alpha < \kappa : \, S \subset A_\alpha \} $ is
countable
and $|A_\alpha \cap A_\beta| \leq \aleph_0$ for $\alpha < \beta <
\kappa$.
Then there are $A_\alpha' \in [A_\alpha]^{\leq \aleph_0}$ for
$\alpha < \kappa$
such that the sequence $\la A_\alpha \bsl A_\alpha'\ra_{\alpha <
\kappa}$
consists of pairwise disjoint sets.
\end{corollary} 

If the existence of certain large cardinals is consistent,
none of the above corollaries is a consequence of GCH alone. To see this, suppose that $\la A_\alpha\ra_{\alpha < 
\kappa}$ is a sequence of sets for which there are $A_\alpha'\in [A_\alpha]^{\aleph_0}$ such that the sequence $\la A_\alpha\setminus A_\alpha'\ra_{\alpha < 
\kappa}$ is point-countable.
 Then on can recursively construct for each $\xi<\kappa$ pairwise disjoint sets $B_\xi\in[\kappa]^{\aleph_0}$ and pairwise disjoint sets $C_\xi\subseteq \bigcup_{\alpha\in B_\xi}A_\alpha\setminus A_\alpha'$ such that $C_\xi\cap A_\alpha\setminus A_\alpha'$ is infinite for each $\alpha\in B_\xi$. Now it is easy to construct a function $c:[\bigcup_{\alpha<\kappa}A_\alpha]^2\rightarrow \{0,1\}$ such that $|c([A_\alpha]^2)|=2$ for all $\alpha<\kappa$. But, 
if the existence of a supercompact cardinal is consistent, 
one can construct a model of ZFC where GCH holds and where there
exists a sequence $\la A_\alpha\ra_{\alpha < 
\aleph_{\omega + 1}}$ of sets of size $\aleph_1$ each
such that every two of these
sets have finite intersection, and for every function 
$c: \bigcup_{\alpha < \aleph{\omega + 1}} A_\alpha \rightarrow \{
0, 1 \}$
there exists $\alpha$ with $|c[A_\alpha]| = 1$ (see Theorem~4.6
of  \cite{HJS}).\\

However, assuming only GCH we can get weaker versions of Theorem~\ref{main}
that are also of interest.  For example:

\begin{theorem}\label{3.1} Suppose that $\sigma$ is 
a regular infinite cardinal and assume that $\lambda^\sigma
=\lambda$ for each cardinal $\lambda$ of cofinality greater than
$\sigma$. Suppose that ${\cal A}=\la A_\alpha\ra_{\alpha <
\kappa}$ is a sequence of
(not necessarily distinct) sets such that 

\noindent $(*)$ For every  $I \in [\kappa]^{\aleph_0} \mbox{ there is a 
finite }
J \subseteq I \mbox{ such that } |\bigcap_{\alpha \in J}
A_\alpha| \leq
\sigma.$

Then there exists $\la A_\alpha'\ra_{\alpha < \kappa}$ such that
$|A_\alpha'| \leq \sigma^+$ for each $\alpha < \kappa$ and
the sequence $\la A_\alpha \bsl A_\alpha'\ra_{\alpha < \kappa}$
is 
point-finite.
\end{theorem}

{\bf Proof:} Assume by induction that $\kappa$ is minimal such that there is 
a counterexample to the theorem. First we show that without loss of generality we may assume that 
$\sigma^+<\kappa$. If not, fix an elementary submodel $M$ of size $\sigma^+\geq \kappa$ such that ${\cal A}\cup\{{\cal A}\}\cup\sigma^+\cup\{\sigma^+\}\subseteq M$.   For each $\alpha\in \kappa$
let $A_{\alpha}'=A_\alpha\cap M$. To see that $\la A_\alpha\setminus A_{\alpha}'\ra_{\alpha < \kappa}$ is point-finite, fix $x\not\in M$ such that $I=\{\alpha:x\in A_{\alpha}\}$ is infinite. By ($*$), fix a finite $J\subseteq I$ such that $|\bigcap_{\alpha \in J}
A_\alpha| \leq
\sigma.$ Clearly $\bigcap_{\alpha \in J}
A_\alpha\in M$ But $\sigma^+\subseteq M$ implies that $\bigcap_{\alpha \in J}
A_\alpha\subseteq M$ contradicting $x\not\in M$.

We may also assume that $\kappa=\lambda^+$ for $\lambda$ a 
cardinal of cofinality $\leq \sigma$. Otherwise the proof is easier (see the comment after the next paragraph). 

Fix a continuous $\in$-chain $\{M_{\xi}:\xi<\kappa\}$ of elementary submodels
of some $H(\Theta)$ such that 
\begin{enumerate}
\item[(a)] $|M_{\xi}|=\lambda$ and $\lambda\subseteq M_{\xi}$ for each $\xi\in\kappa$.
\item[(b)] ${\cal A}\in M_0$ and ${\cal A}\subseteq\bigcup_{\xi\in\kappa}M_{\xi}$.
\item[(c)] For $\xi$ a successor ordinal, $M_{\xi}=\bigcup_{i<cf(\lambda)}M_{\xi}^i$ 
where for each $i$,
\begin{enumerate}
\item[(i)] $\sigma<|M_\xi^i|<\lambda$ and 
\item[(ii)] $[M_\xi^i]^\sigma\subseteq M_\xi^{i+1}$.
\end{enumerate}
\end{enumerate}
Clearly by our assumption such a sequence exists. 

(In the case that $\kappa$ is not a successor of a singular cardinal of 
cofinality $\leq\sigma$ we may fix a sequence as above with the 
stronger property that  $[M_{\xi}]^\sigma\subseteq M_{\xi}$ for each successor $\xi$.)

Notice that as before, by ($*$), we have that for any model $M$ as above and for any 
$x\not\in M$, $\{\alpha\in M:x\in A_{\alpha}\}$ is finite. Indeed, if not, we can
fix a finite subset $J$ for which $|\bigcap_{\alpha\in J}A_{\alpha}|\leq\sigma$.
But then $S=\bigcap_{\alpha\in J}A_{\alpha}$ is an element of  $M$ but it is not a subset. This
implies that $|S|\geq \lambda>\sigma$. A contradiction. 

The main lemma we need
 is the following.

\begin{lemma} $|A_{\alpha}\cap M_\xi|\leq\sigma^+$
for each $\xi<\kappa$ and each $\alpha\not\in M_{\xi}$.
\end{lemma}

{\bf Proof:} Suppose not. First consider the case that $\xi$ is a successor.
Then there is an $i<cf(\lambda)$ such that $|A_{\alpha}\cap M_{\xi}^i|>\sigma^+$. 
Therefore there is an $S\in[A_{\alpha}\cap M_{\xi}^i]^\sigma$ such that
$S\in M_\xi^{i+1}\subset M_{\xi}$. But then $\{\alpha:S\subseteq A_{\alpha}\}
\in M_{\xi}$ but it is not a subset. Therefore it is infinitec contradicting
($*$).
In the case that $\xi$ is a limit consider two subcases: If
$cf(\xi)=\sigma^{++}$, then there is a successor $\eta<\xi$ such that 
$|A_{\alpha}\cap M_{\eta}|\geq\sigma^+$. If 
$cf(\xi)\not =\sigma^{++}$, then there is a successor $\eta<\xi$ such that
$|A_{\alpha}\cap M_{\eta}|\geq\sigma^{++}$. Reasoning as above we get
a contradiction to ($*$) in both cases.\hfill$\Box$\\

Now, to complete the proof of the theorem, using the Lemma and our inductive assumption fix $A_\alpha'\in[A_{\alpha}]^{\leq \sigma^+}$ for each $\xi\in \kappa$ 
and $\alpha\in M_{\xi+1}\setminus M_\xi$
such that 
\begin{enumerate}
\item[(d)] $\{A_\alpha\setminus A_\alpha':\alpha\in M_{\xi+1}\setminus M_\xi\}$ is point-finite, and
\item[(e)] $A_\alpha\cap M_\xi\subseteq A_\alpha'$.
\end{enumerate}
The proof that $\{A_\alpha\setminus A_\alpha':\alpha\in\kappa\}$ is point-finite
is now straightforward: Fix any $x\in\bigcup{\cal A}$. If $x\not\in \bigcup_\xi M_\xi$ then consider
$I=\{\alpha:x\in A_{\alpha}\}$. If $I$ is infinite, fix a finite $J\subseteq I$ such that
$|\bigcap_{\alpha\in J}A_{\alpha}|\leq\sigma$. Then there is a $\xi$ such that $J\subseteq M_\xi$. 
But then $S=\bigcap_{\alpha\in J}A_{\alpha}\in M_\xi$ and $S\subseteq M_{\xi}$ since $|S|\leq \sigma$ contradicting $x\in S\setminus M_{\xi}$.

If $\xi$ is such that $x\in M_{\xi+1}\setminus M_{\xi}$ then
by our construction $\{\alpha:
x\in A_{\alpha}\}\subseteq M_{\xi}$ must be finite as above.
\hfill$\Box$\\

A slightly more general version of Theorem~\ref{3.1} with a slightly weaker conclusion
also follows from GCH alone:

\begin{theorem}\label{3.2} 
Assume GCH, let  $\sigma, \tau$ be
regular infinite cardinals.  Suppose that ${\cal A}=\la A_\alpha\ra_{\alpha <
\kappa}$ is a sequence of
(not necessarily distinct) sets such that 

\noindent ($*$) For every  $I \in [\kappa]^\tau \mbox{ there is }
J \in [I]^{< \tau} \mbox{ such that } |\bigcap_{\alpha \in J}
A_\alpha| \leq \sigma.$

Then there exist $\la A_\alpha'\ra_{\alpha < \kappa}$ such that
$|A_\alpha'| \leq \sigma^+$ for each $\alpha < \kappa$ and
the sequence $\la A_\alpha \bsl A_\alpha'\ra_{\alpha < \kappa}$
is point-$< \tau^+$.
\end{theorem}

\noindent 
{\bf Proof:} We modify the proof of Theorem~\ref{3.1}.  First note that we may assume that $\sigma^+, \tau<\kappa$ and that $\kappa$ is the successor of a limit cardinal $\lambda$ of cofinality $\leq\sigma$. Fix $M_\xi$ satisfying (a), (b) and (c) of the previous proof. In addition, assume

\begin{enumerate}
\item[(f)] $[M^i_\xi]^{\leq \tau}\subseteq M^{i+1}_\xi$, for each successor $\xi<\kappa$ and each $i<cf(\lambda)$.
\end{enumerate}
Define $A_\alpha'$ as before also satisfying $\la A_\alpha \bsl A_\alpha':\alpha\in M_{\xi+1}\setminus M_\xi\ra$ is point-$<\tau^+$. To see that the whole sequence $\la A_\alpha \bsl A_\alpha'\ra_{\alpha < \kappa}$
is point-$< \tau^+$, fix $a\in\bigcup_{\alpha<\kappa} A_\alpha$ and fix $\xi$ maximal so that $a\not\in M_{\xi}$ (where $M_\kappa=\bigcup_{\xi\in \kappa}M_\xi$). It suffices to prove

\begin{lemma} $I=\{\alpha\in M_\xi:a\in A_\alpha\}$ has cardinality $<\tau^+$.
\end{lemma}

{\bf Proof:} Suppose not.\\

{\bf Case 1:} $cf(\lambda)\not=\tau$. Then there are a successor $\eta<\xi$ and an $i<cf(\lambda)$ such that $|I\cap M^i_\eta|\geq \tau$. Fix $J\in [I\cap M_\eta^i]^{<\tau}$ such that $|\bigcap _{\alpha\in J} A_{\alpha}|\leq \sigma$. By (f), $J\in M_\eta$. Therefore $\bigcap _{\alpha\in J} A_{\alpha}\in M_\eta$. Therefore $\bigcap _{\alpha\in J} A_{\alpha}\subseteq M_\eta$ since it is of size $\leq\sigma$. This contradicts $a\not\in M_\xi$. \\

{\bf Case 2:} $cf(\lambda)=\tau$. In the case that $cf(\xi)\not=\tau^+$ we can find a successor $\eta<\xi$ such that $|I\cap M^i_\eta|\geq \tau^+$. This gives a contradiction as in Case 1.
Otherwise assume that $cf(\xi)=\tau^+$. First fix a successor $\eta$ such that $|I \cap M_{\eta}|\geq \tau$. Then fix $J\in[I \cap M_{\eta}]^{<\tau}$ such that $|\bigcap _{\alpha\in J} A_{\alpha}|\leq \sigma$. Now find $i<cf(\lambda)=\tau$ such that $J\subseteq M_\eta^i$. So $J\in M_\eta$ giving a contradiction as before.
\hfill$\Box \Box$\\

Taking $\sigma=\aleph_0$ in Theorem~\ref{3.1} gives the following.

\begin{corollary}\label{GCH}
Assume $\lambda^\omega=\lambda$ for each cardinal $\lambda$ of
uncountable cofinality. Suppose that $\la A_\alpha\ra_{\alpha <
\kappa} $ is a sequence 
of sets such that $|A_\alpha\cap A_{\beta}| \leq \aleph_0$ for each
pair of distinct $\alpha, \beta\in \kappa$. Then
there are $A_\alpha' \in [A_\alpha]^{\leq \aleph_1}$ for $\alpha
< \kappa$
such that the sequence $\la A_\alpha \bsl A_\alpha'\ra_{\alpha <
\kappa}$
is point-finite.
\end{corollary}

Finally, we remark that Theorem~\ref{main} cannot be improved by
demanding that the sets $A_{\alpha }'$ be of cardinality $<\sigma
$. Indeed, if $\sigma$ is regular, then let
$\{A_{\alpha}:\alpha<\sigma^+\}$ be a family of subsets of
$\sigma$ such that $|A_{\alpha}\cap A_{\beta}|<\sigma$ while
$|A_{\alpha}|=\sigma$ for each $\alpha<\sigma^+$. No such family
can be point-$<\sigma^+$. So no point-$<\sigma^+$ family can be
obtained by deleting sets of size less than $\sigma$ from each
$A_{\alpha}$.

\section{Topological Applications}

In this section, we present some topological applications of the
results 
just proved. We are particularly interested in the question,
raised in \cite{HL}, of whether a first countable space with a
weakly uniform
base must have a point-countable base.

\begin{definition}
{\em If $\nio$ and ${\cal B}$ is a base for a topological space $X$,
then we say 
${\cal B}$ is an $n$-{\em weakly uniform base for} $X$ provided
that if 
$A \subseteq X$ with $|A| = n$, then 
$\{ B \in {\cal B} : \; A \subseteq B \}$ is finite. 
We say a base ${\cal B}$ for $X$ is a 
$< \om$-{\em weakly uniform base} provided that if $A \subseteq
X$
with $|A| \geq {\aleph_0}$, then there is a finite subset $F
\subseteq A$
with $\{ B \in {\cal B} :\, F \subseteq B \}$ finite.} 
\end{definition}

The notion of weakly uniform base, introduced in
\cite{HL}, is just what we have called
$2$-weakly uniform base. We shall see later in this section that
the
properties defined above are all distinct.
To avoid trivialities we will be interested only in the case
$n \geq 2$.
Clearly, for $n < m$, $n$-weakly uniform base implies 
$m$-weakly uniform base which implies $< \om$-weakly uniform
base.

It was shown in \cite{HL} that if $X$ is a space with a
$2$-weakly 
uniform base ${\cal B}$ 
and  $x \in X$ is in the closure of a countable subset of 
$X \bsl \{x\}$, then ${\cal B}$ is point-countable at $x$.
In particular, if $X$ is a first countable space with a weakly
uniform 
base~${\cal B}$, then ${\cal B}$ 
is point-countable at all nonisolated points of $X$.
Using that result, it was shown in \cite{HL} that a first
countable space in 
which the boundary of the set of isolated points is separable has
a 
point-countable base.

In \cite{DRW}, it is shown that a first countable space with a
weakly uniform base 
and no more than $\aleph_1$
isolated points has a point-countable base. In \cite{AJRS} this
is improved 
(consistently) to the result that assuming CH, every first
countable 
space with a weakly uniform base and no more than $\aleph_\om$
isolated points
has a point-countable base.

An example is constructed in \cite{DRW}, assuming MA and 
$\ah_2 < 2^{\ah_0}$, of a normal Moore space with a weakly 
uniform base which has no point-countable base.
Such a space also could not be metalindel\"of.

Our next result establishes the independence of the existence of
first 
countable spaces with weakly uniform bases but without 
point-countable bases, and it removes the cardinality restriction
on
the set of isolated points.

\begin{theorem}\label{wun}
Assume CECA. If $X$ is a $T_1$-space in which each non-isolated
point 
is a cluster point of a countable set and $X$ has a $<
\om$-weakly uniform
base, then $X$ has a point-countable base.
\end{theorem}

\noindent
{\bf Proof:}
Suppose ${\cal B}$ is a $< \om$-weakly uniform base for $X$.
Suppose $x$ is a non-isolated point of $X$. We will show that
${\cal B}$
is point-countable at $x$.
Choose a countable set $C \subseteq X \bsl \{ x \}$ with $x \in
\bar{C}$.
For each $B \in {\cal B}$ with $x \in B$, the set $C \cap B$ is
infinite,
and hence there exists $F_B \in [ C \cap B ]^{<\ah_0}$
such that $\{ A \in B \; : F_B \subset A \}$ is finite.
Since $[ C ]^{< \ah_0}$ is countable and for each
$F \in [ C ]^{< \ah_0}$ there are only finitely many 
$B \in {\cal B}$ with $F_B = F$, 
we conclude that 
$\{ B \in {\cal B}  \; : x \in B \}$ is countable.

Let $\la x_\alpha \; : \al < \kappa \ra$ list the isolated points
of 
$X$, and let $A_\al = \{ B \in {\cal B} \; : x_\al \in B \}$. 
Notice that since ${\cal B}$ is a 
$< \om$-weakly uniform base, we have that if $I$ is an
infinite subset of $\kappa$, then there exists a finite $J
\subseteq I$
such that 
$| \{ B \in {\cal B} \; : \{ x_\al \; : a \in J \} \subseteq B \}
| < {\aleph_0}$,
i.e.
$| \bigcap_{\al \in J} A_\al | < {\aleph_0}$.
By Corollary~\ref{cor1}, 
for each $\al$ there is $A^1_\al \in [ A_\al ]^{\leq {\aleph_0}}$
such that 
$\la A_\al \bsl A_\al^1 \; : \al < \kappa \ra$ is point-finite,
i.e. for 
each $B \in {\cal B}$ we have 
$B \in A_\al \bsl A_\al^1$ for only finitely many $\al$'s. 
Hence for each $B \in {\cal B}$, there 
is a finite set $I(B)$ of isolated points
such that for every $\al < \kappa$, 
$\{ B \in {\cal B} \; : x_\al \in B \bsl I(B) \} \subseteq
A^1_\al$,
and $| A_\al^1 | \leq {\aleph_0}$.
Hence ${\cal B} = \{ B \bsl I(B) \; : B \in {\cal B} \}$ is a 
point-countable open family of subsets of $X$ which forms a base
at every
non-isolated point of $X$. Thus 
${\cal B}^* = {\cal B}_1 \cup \{ \{ x_\al \} \; : \al < \kappa \}$ is
a point-countable base for $X$.
\hspace{\fill} $\Box$

\begin{corollary}
It is consistent with $ZFC$ that every first countable
$T_1$-space
with a weakly uniform base has a point-countable base.
\end{corollary}

Combining that corollary with the example of \cite{DRW} completes
the 
independence result.

Assuming only $\lambda^{\omega }=\lambda $ for each cardinal
$\lambda $ of uncountable cofinality, we can obtain a weaker
version of Theorem~\ref{wun}. The proof is identical, using
Corollary~\ref{GCH} in place of Corollary~\ref{cor1}.

\begin{theorem}
Assume $\lambda^{\omega  }=\lambda $ for each cardinal $\lambda $
of uncountable cofinality. If $X$ is a $T_1$-space in which each
non-isolated point 
is a cluster point of a countable set and $X$ has a $<
\om$-weakly uniform
base, then $X$ has a point-$<\aleph _{2}$ base.
\end{theorem}

\begin{corollary} Assume $\lambda^{\omega  }=\lambda $ for each cardinal $\lambda $
of uncountable cofinality. Every first countable
$T_1$-space
with a weakly uniform base has a point-$<\aleph_2$ base.
\end{corollary}

\medskip

\begin{definition}
{\em If $\nio$, then we say that a topological space $X$ is {\em
n-metacompact}
provided that for every open cover ${\cal U}$ 
of $X$ there is an open refinement 
${\cal V} \prec {\cal U}$ 
such that if $A \subseteq X$ with $| A | = n$, then 
$\{ V \in {\cal V} \; : A \subseteq V \}$ is finite.
We say $X$ is $< \om$-metacompact provided that for every open
cover ${\cal U}$ of $X$ there is an open refinement ${\cal V}
\prec {\cal U}$
such that if $A \subseteq X$ with $| A | \geq {\aleph_0}$, then
there is a finite subset $F \subseteq A$ with 
$\{ V \in {\cal V} \; : F \subseteq V \}$ finite.}
\end{definition}

It is clear that $1$-metacompact is just metacompact and that for
$n <m$, $n$-metacompact implies $m$-metacompact which implies 
$< \om$-metacompact. We will call a refinement
$\cal{V}$ as in the definition above an {\em n-weakly uniform
refinement},
or {\em $< \om$-weakly uniform refinement}
respectively.

\begin{theorem}
Assume CECA. If $X$ is a $< \om$- metacompact $T_1$-space, then
$X$ is metalindel\"of.
\end{theorem}

\noindent
{\bf Proof:}
Suppose $X$ is $< \om$-metacompact and $\cal{U}$ is an open cover
of 
$X$. Let ${\cal V}$ be an open, $< \om$-weakly uniform refinement
of $\cal{U}$.
List the points of $X$ as 
$\la x_\al \; : \al < \kappa \ra$.
For each $\al <\kappa$, let 
$A_\al = \{ V \in{\cal V}\; : x_\al \in V \}$. 
By Corollary~\ref{cor1}, for each $\al <\kappa$ there is 
$A_\al^1 \in [A_\al]^{\leq {\aleph_0}}$ such that 
$\la A_\al \bsl A_\al^1 \; : \al \in\kappa \ra$ is point-finite
on the 
set ${\cal V}$. For each $V \in{\cal V}$, we let 
$I(V) = \{ x_\al \in V \; : V \in A_\al \bsl A_\al^1 \}$. 
Note that 
$| I(V) | < {\aleph_0}$.
Further, for each $\al <\kappa$, 
$\{ V \in {\cal V}\; : x_\al \in V \bsl I(V) \} \subseteq
A^1_\al$ and
$| A_\al^1 | \leq {\aleph_0}$. For each $\al <\kappa$, choose 
$V_\al \in {\cal V}$ such that
$x_\al \in V_\al$. For $V \in{\cal V}$, let 
$V^* = ( V \bsl I(V)) \cup \{ x_\al \in I(V) \; : V = V_\al \}$.
Note that $V^* \subseteq V$ and $V \bsl V^*$ is finite. Let
${\cal V}^* = \{ V^* \; : V \in {\cal V}\}$.
Then ${\cal V}^*$ is an open refinement of $\cal{U}$ and if 
$x_\al \in V^*$, then 
$V \in A_\al^1 \cup \{ V_\al \}$.  Thus ${\cal V}^*$ is
point-countable.

\begin{corollary}
It is consistent with $ZFC$ that every $T_1$-space with a weakly 
uniform base is metalindel\"of.
\end{corollary}

Again, combining this corollary with the example in \cite{DRW}
completes the 
independence result.

\begin{example}\label{ex1}
For each natural number $n \geq 1$, there is a Moore space $X_n$ 
of scattered height 2 which has an 
$(n + 1)$-weakly uniform base, but $X_n$ does not have an
$n$-weakly uniform base.
\end{example}

\noindent
{\bf Proof:}
Let $L$ be a subset of $\bR \times \{ 0 \}$ with $|L |=
{\aleph_1}$, 
and let $D = \{p_n \; :n \in \om \}$ be a countable subset of 
$\bR \times ( 0, \infty)$ which is dense in the Euclidean
topology.
It is shown in \cite{DRW} that there is a collection $\cal{H}$ of
countably 
infinite subsets of $L$ and a partition 
$\{ {\cal H}_n \; : \nio \}$ of ${\cal H}$ such that 
\begin{enumerate}
\item[(1)]
if $H_1$, $H_2 \in \cal{H}$ and $H_1 \neq H_2$, then 
$| H_1 \cap H_2 | < {\aleph_0}$, and
\item[(2)]
if $Y \subseteq L$ and $|Y |= \aleph_1$, 
then for each $\nio$, there exists $H \in {\cal H}_n$ such that 
$| Y \cap H | = {\aleph_0}$.
\end{enumerate}
For each $\nio$, let
$K_n = \{ ( p_n , H ) \; : H \in {\cal H}_n \}$ and let 
$K = \bigcup_{\nio} K_n$.
Let $X = L \cup K$. For each $\nio$ and each $x \in L$, 
let $B_n (x) = \{ (p_i , \, H ) \; : ( p_i , \, H ) \in K_i$, $x
\in H$, 
and $p_i$ is an element of the Euclidean open ball in 
$\bR  \times (0, \infty )$ of radius 
$2^{-n}$ which is tangent to the axis 
$\bR \times \{ 0 \}$ at the point $x \}$.
If $x \in L$ and $\nio$, then we define $G_n (x) = \{ x \} \cup
B_n (x)$
and let 
$\{ G_n (x) \; : \nio \}$ be a neighborhood base at $x$. 
If $y \in K$, then $\{ y \}$ is open. It 
is shown in \cite{DRW} that with this topology $X$ is a Moore
space 
with a weakly uniform base, and if $\cal{U}$ is any open cover of
$X$
which refines the open cover
$\{ G_0 (x) \; : x \in L \} \cup \{ \{ y \} \; : y \in K \}$,
then
there exists $y \in K$
such that $\{ U \in {\cal U} \; : y \in U \}$ is infinite.

Suppose $k$ is a natural number and $k \geq 1$. 
We shall use the space $X$ to construct a Moore space $X_k$ of 
scattered height $2$ with a $(k+1)$-weakly uniform base but no
$k$-weakly
uniform base.
Let $X_k = L \cup (K \times k)$.
Points of $K \times k$ will be isolated. If $x \in L$ and $\nio$,
let $G^k_n (x) = \{ x \} \cup ( B_n (x) \times k )$. 
Notice that $\la G_n^k (x) \; : \nio \ra$ 
is a decreasing sequence, and thus this is a valid assignment 
of neighborhoods. 
Letting ${\cal G}_n = \{ G_n^k (x)\; : x \in L \} \cup \{ \{ y
\} : y \in K \times k \}$, 
we see that $\la {\cal G}_n \; : \nio \ra$ is a development for
$X_k$, and 
since each $G_n^k (x)$ is clopen, 
$X_k$ is a $0$-dimensional space of scattered height $2$.

We now show that ${\cal B} = \bigcup_{\nio} {\cal G}_n$ is a 
$(k+1)$-weakly uniform base for $X_k$.
Suppose $A \subseteq X_k$ with $|A |= k+1$. 
If $| A \cap L | \geq 2$, then 
$|\{ B \in {\cal B}  : A \subseteq B \} |= 0$. 
If $| A \cap L | = 1$. Then let $\{ x \} = A \cap L$.
Choose $\nio$ such that $G_n^k (x) \cap A = \{ x \}$.
Then $|\{ B \in {\cal B} :\, A \subseteq B \} |\leq n$,
since for each $m \in \om$ there is only one element of ${\cal
G}_m$ which
contains $x$. 
Finally, if $|A \cap L |=0$, then let
$A = \{ (y_1 , n_1 ) , \dots , ( y_{k+1} , n_{k+1}) \}$
where $y_i \in K$ and $n_i \in k$ for $1 \leq i \leq k+1$.
Choose $i$, $j \leq k+1$ such that $n_i = n_j$.
Then $\{ B \in {\cal B} \; : A \subseteq B \} \subseteq \{ B \in
{\cal B} \; :
\{ (y_i , n_i ) , ( y_j , n_j) \} \subseteq B \}$.
Now for $B = G^k_n (x)$, we have that 
$( y_j , n_j) \in  G^k_n (x)$ implies $y_i \in G_n (x)$.
So 
$|\{ B \in {\cal B} \; : \{ (y_i , n_i ) , ( y_j , n_j) \}
\subseteq B \} |
= | \{ G_n (x) \; : \{ y_i , y_j  \} \subseteq G_n (x) \} |<
{\aleph_0}$,
since $\{ G_n (x) \; : \nio , x \in L \}$ is shown in \cite{DRW}
to be a 
weakly uniform collection.

We now show that $X_k$ has no $k$-weakly uniform base.
Suppose ${\cal V}$ 
is any open cover of $X_k$ which refines ${\cal
G}_1$.
For each $x \in L$, choose $V_x \in{\cal V}$ and $n (x) \in \om$
with 
$G^k_{n(x)} (x) \subseteq V_x$.
Now in the space $X$, ${\cal U} = \{ G_{n(x)} (x) \; : x \in L \}
\cup \{ \{ y \} \; : y \in K \}$
refines 
$\{ G_1 (x) \; : x \in L \} \cup \{ \{ y \} \; : y \in K \}$.
Hence there is a point
$y \in K$ so that $\{ U \in {\cal U} \; : y \in U \}$ is
infinite, and so 
$\{ x  \in L \; : y \in G_{n(x)} (x) \}$ is infinite.
Let $A = \{ y \} \times k$. 
 Then $|A |= k$
and for each $i \in k$, $(y, \, i ) \in G_{n(x)}^k (x) \subseteq
V_x$.
So $\{ x \in L \; : A \subseteq V_x \}$ 
is infinite, and thus ${\cal V}$ is not a $k$-weakly uniform base
for 
$X_k$.
\hspace{\fill} $\Box$

\begin{example}
There is a space $Y$ which has $< \om$-weakly uniform base, but
does not 
have $n$-weakly uniform base for any $\nio$.
\end{example}

\noindent
{\bf Proof:}
For each natural number $n \geq 1$, let $X_n$ be as constructed
in 
Example~\ref{ex1}.
Let $Y$ be the disjoint union of the spaces $X_n$. The natural
base is 
easily seen to be $< \om$-weakly uniform since if $A \subseteq Y$
and 
$|A |= \aleph_0$, then either $A \subseteq X_n$ for some $n$ in
which case
any subset of $A$ of size $n+1$ would be contained in only
finitely many 
elements of the base, or $A$ contains two points from distinct 
$X_n$'s in which case no element of the base contains that
two-point subset.
To see that $Y$ cannot have an $n$-weakly uniform base, it is
enough to observe that $X_n$ is an open subspace of $Y$.
\hspace{\fill} $\Box$\\

Thus we have shown that these properties are all distinct.
We leave the reader with one question regarding covering
properties.

\begin{question}
Is it true or consistent that every space with a weakly uniform
base
must be submetacompact?
\end{question}

\noindent
AUTHOR'S ADDRESSES:\\

Balogh and Davis

Department of Mathematics

Miami University

Oxford, OH 45056

U.S.A.\\

Just and Szeptycki

Department of Mathematics

Ohio University

Athens, OH 45701

U.S.A.\\

Shelah

Institute of Mathematics

The Hebrew University

Givat Ram

91904 Jerusalem

ISRAEL

\end{document}